\documentclass[11pt]{amsart}

\usepackage{amsmath,amssymb,amsthm,mathtools}
\usepackage[margin=1.1in]{geometry}
\usepackage{enumitem,placeins,tikz}
\usetikzlibrary{positioning}
\tikzset{latticepoint/.style={circle,fill=black,inner sep=1.5pt}}
\usepackage[hidelinks]{hyperref}

\usepackage{orcidlink}

\theoremstyle{plain}
\newtheorem{theorem}{Theorem}[section]
\newtheorem{lemma}[theorem]{Lemma}
\newtheorem{proposition}[theorem]{Proposition}
\newtheorem{corollary}[theorem]{Corollary}
\newtheorem*{theoremA}{Theorem A}
\newtheorem{conjecture}{Conjecture}
\theoremstyle{definition}
\newtheorem{definition}[theorem]{Definition}
\theoremstyle{remark}

\newcommand{\Z}{\mathbb{Z}}
\newcommand{\R}{\mathbb{R}}
\newcommand{\C}{\mathbb{C}}

\newcommand{\one}{\mathbf{1}}

\newcommand{\supp}{\operatorname{supp}}
\newcommand{\ellone}{\ell^1(\Z^2)}
\newcommand{\Tstar}{\mathcal{T}^{\star}_{\theta}}
\newcommand{\Sfour}{S_4}

\begin{document}

\title[Discrete X-ray transform under quadratic decay]{Counterexamples to uniqueness for a discrete X-ray transform under quadratic decay}

\author[S. R. Jathar]{Shubham R. Jathar\,\orcidlink{0000-0002-9109-7543}}
\address{Computational Engineering, School of Engineering Sciences,
Lappeenranta-Lahti University of Technology LUT, Lappeenranta, Finland}
\email{shubham.jathar7@gmail.com, shubham.jathar@lut.fi}

\author[J. Railo]{Jesse Railo\,\orcidlink{0000-0001-9226-4190}}
\address{Computational Engineering, School of Engineering Sciences,
Lappeenranta-Lahti University of Technology LUT, Lappeenranta, Finland}
\email{jesse.railo@lut.fi}

\date{}
\dedicatory{Dedicated to Roman Novikov on the occasion of his 60th birthday}

\begin{abstract}
We show nonuniqueness for a discrete X-ray transform on $\Z^2$ under quadratic decay. We construct an explicit nonzero function in $\ell^1(\Z^2)$ with pointwise decay $O(|z|^{-2})$ whose sums over integer lines vanish. We also obtain an infinite-dimensional null space with the same decay assumption. Furthermore, we show a result which indicates that the decay rate $|z|^{-2}$ may be sharp and faster decay may imply uniqueness.
\end{abstract}
\maketitle

\section{Introduction}\label{sec:intro}

For a function $u\colon\R^2\to\C$ whose line integrals converge absolutely, the X-ray transform assigns to each straight line $L\subset\R^2$ the integral of $u$ over $L$. These line integrals provide a mathematical model for the data in computerized tomography \cite{Helgason:Book,Sharafutdinov:book,GIP2D}. A discrete analogue replaces each line integral by the sum of the values of a function $f\colon\Z^2\to\C$ at the lattice points on that line. Strichartz \cite[Section~5, Problem~B, pp.~383--384]{Strichartz:1982} defines this discrete Radon transform for absolutely summable functions, using all straight lines that contain at least two lattice points. For a smaller family of measured lines, we ask which support or decay assumptions ensure that the remaining line sums still determine $f$ uniquely.

We study this question for the planar line family introduced by Novikov and Sharma \cite[Section~2.2, (23)--(26)]{Novikov:Sharma:2026}. Their construction associates a line $\gamma_z$ with each lattice point $z\in\Z^2$. We write $|z|$ for the Euclidean norm and $z\cdot y$ for the dot product on $\R^2$. Let $J(z_1,z_2)=(-z_2,z_1)$ denote the counterclockwise rotation through $90^\circ$ about the origin. For $z\in\Z^2\setminus\{0\}$, we define
\begin{equation}\label{eq:gammaz}
\gamma_z:=\{\,y\in\R^2:z\cdot y=|z|^2\,\}=z+\R Jz.
\end{equation}
Thus $\gamma_z$ passes through $z$ and is perpendicular to the vector $z$.

At $z=0$, the expression $z+\R Jz$ reduces to the single point $\{0\}$. We therefore choose the line through the origin separately. Fix a primitive vector $\theta\in\Z^2\setminus\{0\}$ and define \begin{equation}\label{eq:gamma0}
\gamma_0:=\R\theta.
\end{equation}
Here \emph{primitive} means that the greatest common divisor of $|\theta_1|$ and $|\theta_2|$ is $1$. With this choice, the family of measured lines is
\begin{equation}\label{eq:Tstar}
\Tstar:=\{\,\gamma_z:z\in\Z^2\,\}.
\end{equation}
We call $\gamma_z$ the \emph{selected line} at $z$, with \emph{base point} $z$. The family $\Tstar$ depends on $\theta$ only through the line $\gamma_0$.

For a line $L\subset\R^2$ and a function $f\colon\Z^2\to\C$, we define the \emph{line sum}
\begin{equation}\label{eq:linesum}
Pf(L):=\sum_{y\in L\cap\Z^2}f(y),
\end{equation}
provided that the series converges absolutely. For functions whose sums converge absolutely on every selected line, the measured data are the collection
$(Pf(\gamma_z))_{z\in\Z^2}$. 

Novikov and Sharma prove that these data determine $f$ uniquely within the class of finitely supported functions \cite[Theorem~2.1]{Novikov:Sharma:2026}. For a function $f$ supported in a known disk, the reconstruction of Novikov and Sharma recovers $f(z)$ in decreasing order of $|z|$. The known zero values outside the disk provide the starting step. For a function with infinite support, there is no outermost support point from which to begin the same procedure. We ask whether the line sums $Pf(\gamma_z)$, $z\in\Z^2$, uniquely determine $f$ under a suitable decay assumption at infinity, without requiring finite support.

Novikov and Sharma formulate two conjectures related to global and exterior uniqueness under decay assumptions on $f$. For $\varepsilon>0$ and $\alpha\ge0$, we define the space $X_\varepsilon$ and the set $E_\alpha$ by
\begin{equation}\label{eq:Xeps}
\begin{aligned}
X_\varepsilon
&:=\Bigl\{f\colon\Z^2\to\C:
\sup_{z\in\Z^2}(1+|z|)^{1+\varepsilon}|f(z)|<\infty\Bigr\},\\
E_\alpha
&:=\{z\in\Z^2:|z|\ge\alpha\}.
\end{aligned}
\end{equation}
Thus $f\in X_\varepsilon$ means that $|f(z)|\le C_f(1+|z|)^{-1-\varepsilon}$ for every $z\in\Z^2$ and some finite constant $C_f\ge0$.
The constant may depend on $f$ and $\varepsilon$, but not on $z$. The set $E_\alpha$ consists of the lattice points outside the open disk of radius $\alpha$ centered at the origin.

A function $f\colon\Z^2\to\C$ is called \emph{rapidly decreasing} if, for every $N\ge0$,
\[
f(z)=O((1+|z|)^{-N})
\qquad\text{as }|z|\to\infty.
\]
The selected line sums converge absolutely for functions in $X_\varepsilon$ and for rapidly decreasing functions, as proved in Lemma~\ref{lem:convergence}. With these definitions, the conjectures can be stated as follows.
\begin{conjecture}[{\cite[Conjecture~2.3]{Novikov:Sharma:2026}}]
Let $\varepsilon>0$ and let $f\in X_\varepsilon$.
If $Pf(\gamma)=0$ for every $\gamma\in\Tstar$,
then $f=0$.
\end{conjecture}

\begin{conjecture}[{\cite[Conjecture~2.4]{Novikov:Sharma:2026}}]
Let $\alpha\ge0$ and let $f\colon\Z^2\to\C$
be rapidly decreasing.
If $Pf(\gamma_z)=0$ for every $z\in E_\alpha$,
then $f=0$ on $E_\alpha$.
\end{conjecture}
The second conjecture asks whether the line sums with base points in $E_\alpha$ determine $f$ on $E_\alpha$. No values of $f$ inside the open disk are assumed known.

\subsection{Main results}\label{sec:main}
Our counterexamples are supported on the coordinate axes and the lines $z_2=\pm z_1$. We write
\begin{equation}\label{eq:S4}
\Sfour:=\{\,z\in\Z^2: z_1z_2(z_1-z_2)(z_1+z_2)=0\,\}
\end{equation}
for the set of lattice points on these four lines. 

For a function $f\colon\Z^2\to\C$, we define its support by
\[
\supp f:=\{z\in\Z^2:f(z)\ne0\}.
\] We write $\ellone$ for the space of absolutely summable functions on $\Z^2$, equipped with the norm
\[
\|f\|_{\ell^1}:=\sum_{z\in\Z^2}|f(z)|<\infty.
\]
\begin{theoremA}
Fix a primitive vector $\theta\in\Z^2\setminus\{0\}$,
and let $\Tstar$ be the family defined in \eqref{eq:Tstar}.
\begin{enumerate}[label=\textup{(\roman*)},leftmargin=2.2em]
\item\label{A:a}
There exists a nonzero real-valued function
$f_\theta\in\ellone$ with $\supp f_\theta\subseteq\Sfour$
such that
\[
|f_\theta(z)|\le|z|^{-2}\quad(z\ne0),
\qquad
Pf_\theta(\gamma)=0\quad(\gamma\in\Tstar).
\]
For every $\varepsilon>0$, we have
$f_\theta\in X_\varepsilon$ if and only if
$\varepsilon\le1$.

\item\label{A:b}
There exists a linearly independent sequence
$(G_j)_{j\ge0}$ of real-valued functions in
$X_1\cap\ellone$.
For every integer $j\ge0$, we have
\[
G_j(0)=0,
\qquad
PG_j(\gamma_z)=0
\quad(z\in\Z^2\setminus\{0\}).
\]
Moreover, $PG_j(L)=0$ for every $j\ge0$
and every line $L\subset\R^2$ through the origin.

\item\label{A:c}
Let $\alpha\ge0$ and let $f\colon\Z^2\to\C$
satisfy $\supp f\subseteq\Sfour$.
Assume that $f(z)=o(|z|^{-2})$ as $|z|\to\infty$.
If $Pf(\gamma_z)=0$ for every $z\in E_\alpha$,
then $f=0$ on $E_\alpha$.
\end{enumerate}
\end{theoremA}
Part~\ref{A:a} disproves Conjecture~2.3 of \cite{Novikov:Sharma:2026}
for $0<\varepsilon\le1$. The functions in part~\ref{A:b} span an
infinite-dimensional subspace of $X_1\cap\ellone$ whose elements have zero selected line sums. Their sums also vanish on every line through
the origin.

Part~\ref{A:c} proves Conjecture~2.4 of \cite{Novikov:Sharma:2026} for functions supported on $\Sfour$, under the weaker decay assumption
$f(z)=o(|z|^{-2})$. This condition cannot be weakened to $O(|z|^{-2})$. Indeed, the counterexample in part~\ref{A:a} has zero selected line sums and is nonzero arbitrarily far from the origin. These results do not settle Conjecture~2.4 for functions with unrestricted support.

\subsection{Revised Novikov's conjecture}
\label{sec:conjecture}
We ask whether the exterior uniqueness conclusion
of Theorem~A\ref{A:c} holds without the support condition
$\supp f\subseteq\Sfour$.
\begin{conjecture}[Revised Novikov's conjecture]
Let $\alpha\ge0$ and let $f\colon\Z^2\to\C$ satisfy
\[
f(z)=o(|z|^{-2})
\qquad\text{as }|z|\to\infty.
\]
If $Pf(\gamma_z)=0$ for every $z\in E_\alpha$,
then $f=0$ on $E_\alpha$.
\end{conjecture}
All selected line sums in the conjecture converge absolutely by Lemma~\ref{lem:convergence}, with $s=2$. Every rapidly decreasing function satisfies $f(z)=o(|z|^{-2})$. The conjecture would therefore imply Conjecture~2.4 of \cite{Novikov:Sharma:2026}. For $\varepsilon>1$ and $f\in X_\varepsilon$, definition \eqref{eq:Xeps} gives
\[
|z|^2|f(z)|
=O(|z|^{1-\varepsilon})\longrightarrow0
\qquad\text{as }|z|\to\infty.
\]
Thus the case $\alpha=0$ would also establish
Conjecture~2.3 of \cite{Novikov:Sharma:2026}
for $\varepsilon>1$.

\subsection{Related literature}
Decay assumptions also appear in uniqueness results for the continuous X-ray transform. For $n\ge2$, a continuous function
$u\colon\R^n\to\C$ satisfying
\[
|u(x)|\le C(1+|x|)^{-\eta}
\qquad(x\in\R^n)
\]
for some constants $C>0$ and $\eta>1$ is uniquely determined by its integrals over all straight lines.
This follows from Jensen's inversion formula \cite[Theorem~5.3]{Jensen:2004} with $k=1$ and $\delta=\eta-1$. 
For classical injectivity results, see \cite{Helgason:Book}. Related results for geodesic X-ray transforms on Cartan--Hadamard manifolds are proved in \cite{Lehtonen:Railo:Salo:2018}.

Decay along each line separately does not ensure uniqueness. Zalcman \cite[Section~5, pp.~243--244]{Zalcman:1982} constructed a nonzero continuous function on $\R^2$ whose integrals over all straight lines converge absolutely and vanish. The function is $O(|x|^{-2})$ along each fixed line, but this decay estimate does not hold uniformly on $\R^2$. For a related nonuniqueness example, see \cite{Armitage:1994}. Our counterexample satisfies the uniform bound $|f_\theta(z)|\le|z|^{-2}$ for every nonzero $z\in\Z^2$, but only the selected line sums over $\Tstar$ are measured.

The choice of weight also affects injectivity. Boman \cite{Boman:1993} constructed a smooth strictly positive weight depending on position
and line direction. For this weight, a nonzero smooth compactly
supported function on $\R^2$ has zero weighted integrals over all straight lines. Goncharov and Novikov \cite{Goncharov:Novikov:2019}
constructed corresponding counterexamples in higher dimensions with strictly positive continuous weights. 

Nonuniqueness also occurs in Fourier phase retrieval. Novikov and Xu \cite{Novikov:Xu:2026} constructed nonuniquness examples for recovery from Fourier magnitude in discrete and continuous settings. Their measurement are nonlinear, whereas our selected line sum are linear.

Strichartz \cite[Section~5, Problem~D, p.~384]{Strichartz:1982}
posed a related reconstruction problem on the flat torus $\mathbb T^2=\R^2/\Z^2$, using mean values over closed geodesics. Abouelaz and Rouvi\`ere \cite{Abouelaz:2011:tori} obtained an inversion formula for continuous functions with absolutely summable Fourier
coefficients. Ilmavirta \cite[Theorem~2]{Ilmavirta:2015}
proved injectivity of periodic Radon transforms on distributions. Railo \cite{Railo:2020} developed reconstruction formulas and stability estimates for periodic $d$-plane transforms on $\mathbb T^n$, $1\le d<n$. For numerical reconstruction, Ilmavirta, Koskela, and Railo \cite{Ilmavirta:Koskela:Railo:2020}
developed Torus CT, which Salo et al. \cite{Salo:Meaney:Koskela:Railo:2026} adapted to experimental X-ray data. These torus transforms use integrals over closed geodesics or periodic planes, whereas our measurements sum function values over lattice
points on the selected lines $\gamma_z$.

Other discrete formulations use different measurement operators. Fokas and Novikov \cite{Fokas:Novikov:1991} introduced a discrete analogue of the Radon transform based on finite-difference equations.
The approximate discrete Radon transform (ADRT) sums pixel values along digital lines that approximate straight lines
\cite{Brady:1998,Press:2006}. Exact inversion and range characterizations
for the ADRT were obtained in \cite{Rim:2020,Li:Ren:Rim:2023,Rim:2026}.

Gupta and Novikov \cite{Gupta:Novikov:2026} study a discrete non-abelian X-ray transform defined by ordered matrix products. For the same planar selected line family, they prove uniqueness for invertible matrix-valued
functions equal to the identity outside a finite set and pose related questions under decay towards the identity
\cite[Theorem~2.1 and Questions~1--2]{Gupta:Novikov:2026}.

It would be interesting to study analogous decay questions for lattice sums over restricted families of rational affine $d$-planes in $\R^n$.
\section{Preliminaries}\label{sec:prelim}
We describe the lattice points on the selected lines and establish absolute convergence of their line sums.
\subsection{Vectors and rotation}
We write points of $\R^2$ as $z=(z_1,z_2)$. The lattice is
$\Z^2=\{z\in\R^2:z_1,z_2\in\Z\}$. Every disk of radius $R\ge0$ centered at the origin
contains finitely many lattice points. Indeed, $|z|\le R$ implies $|z_i|\le R$ for $i=1,2$, so each integer coordinate has only finitely many possible values.
We write $I$ for the identity map on $\R^2$. We will use the following properties of the rotation $J$.
\begin{lemma}\label{lem:J}
For all $z,y\in\R^2$, the following hold.
\begin{enumerate}[label=\textup{(\roman*)},leftmargin=2.2em]
\item
$J^2=-I$, $J^4=I$, and $J^{-1}=-J=J^3$.
\item
$Jz\cdot Jy=z\cdot y$, $|Jz|=|z|$,
and $z\cdot Jz=0$.
\item
$J$ maps $\Z^2$ bijectively onto $\Z^2$.
\end{enumerate}
\end{lemma}

\begin{proof}
Direct calculation gives
\[
J^2z=-z,\qquad
Jz\cdot Jy=z_2y_2+z_1y_1=z\cdot y,
\qquad
z\cdot Jz=0.
\]
The first identity implies $J^4=I$ and $J^{-1}=-J=J^3$. Taking $y=z$ in the dot product identity gives $|Jz|=|z|$. Both $J$ and its inverse $-J$ send integer coordinates to integer coordinates. Hence $J$ maps $\Z^2$ bijectively onto $\Z^2$.
\end{proof}
\subsection{Primitive integer vectors}
For an integer $a$ and a positive integer $n$, we write $n\mid a$ if $a=nq$ for some $q\in\Z$.
\begin{definition}\label{def:primitive}
For $z\in\Z^2\setminus\{0\}$, let
\[
g(z):=\gcd(|z_1|,|z_2|)
\]
be the greatest positive integer dividing both coordinates. If one coordinate is zero, then $g(z)$ is the absolute value of the other coordinate. The vector $z$ is called \emph{primitive} if $g(z)=1$.
\end{definition}
We use primitive integer vectors to parametrize the lattice points on the selected lines. The following lemma records the arithmetic facts needed for this parametrization.
\begin{lemma}\label{lem:bezout}
Let $z\in\Z^2\setminus\{0\}$, and let $a\in\Z^2$ be primitive.
\begin{enumerate}[label=\textup{(\roman*)},leftmargin=2.2em]
\item\label{bz:a}
The vector $z/g(z)$ has integer coordinates and is primitive.

\item\label{bz:b}
There exist integers $u,v$ such that $ua_1+va_2=1$.

\item\label{bz:c}
If $t\in\R$ and $ta\in\Z^2$, then $t\in\Z$.

\item\label{bz:d}
We have $g(Jz)=g(z)$. In particular, $Ja$ is primitive.
\end{enumerate}
\end{lemma}

\begin{proof}
For part~\ref{bz:a}, put $b:=z/g(z)$. Since $g(z)$ divides both coordinates of $z$,
the vector $b$ has integer coordinates. Moreover, $b\ne0$ because $z\ne0$.
Let $d$ be a positive integer dividing both coordinates of $b$. Then $d\,g(z)$ divides both coordinates of $z$. Since $g(z)$ is their greatest positive common
divisor, we have $d\,g(z)\le g(z)$. Dividing by $g(z)>0$ gives $d\le1$. Since $d$ is a positive integer, $d=1$. Thus $b$ is primitive.

For part~\ref{bz:b}, primitivity of $a$ gives $\gcd(a_1,a_2)=1$. B\'ezout's identity
\cite[Theorem~1.7]{Jones:Jones:book} therefore gives integers $u,v$ such that
$ua_1+va_2=1$.

For part~\ref{bz:c}, choose integers $u,v$ as in part~\ref{bz:b}. Since $ta\in\Z^2$, both $ta_1$ and $ta_2$ are integers. Consequently, we have
\[
t=t(ua_1+va_2)=u(ta_1)+v(ta_2)\in\Z.
\]

For part~\ref{bz:d}, the coordinates of $Jz=(-z_2,z_1)$ differ from those of $z$
only by their order and signs. Hence the positive integers dividing both coordinates
are the same for $z$ and $Jz$, which gives $g(Jz)=g(z)$. Since $a$ is primitive, we have $g(Ja)=g(a)=1$. Thus $Ja$ is primitive.
\end{proof}

\subsection{Lines through lattice points}
Let $p\in \Z^2$ and $d\in\Z^2\setminus\{0\}$. The vector $w:=d/g(d)$ is primitive and satisfies $p+\R d=p+\R w$. The following lemma shows that the lattice points on this line are exactly the points $p+kw$,
with $k\in\Z$.
\begin{lemma}\label{lem:rational}
 Let $p\in\Z^2$ and let $w\in\Z^2$ be primitive. Then
\[
(p+\R w)\cap\Z^2=\{\,p+kw : k\in\Z\,\},
\]
and distinct integers $k$ give distinct points.
\end{lemma}

\begin{proof}
Let $y\in(p+\R w)\cap\Z^2$. Then $y=p+tw$ for some $t\in\R$, and $tw=y-p\in\Z^2$. Since $w$ is primitive, Lemma~\ref{lem:bezout}\ref{bz:c} gives $t\in\Z$. Conversely, for every $k\in\Z$, the point $p+kw$ has integer coordinates and lies on $p+\R w$.
Thus the lattice points on $p+\R w$ are exactly the points $p+kw$ with $k\in\Z$. If $p+k_1w=p+k_2w$ for integers $k_1,k_2$, then $(k_1-k_2)w=0$. Since $w\ne0$, we obtain $k_1=k_2$. Thus distinct integers give distinct lattice points.
\end{proof}
\subsection{The selected family of lines}\label{sec:selected}
Recall the selected lines $\gamma_z$ and the family $\Tstar$ defined in \eqref{eq:gammaz}--\eqref{eq:Tstar}. To parametrize their lattice points, we define the direction vectors
\begin{equation}\label{eq:vz}
v_z:=\frac{Jz}{g(z)}
\quad(z\in\Z^2\setminus\{0\}),
\qquad
v_0:=\theta.
\end{equation}
The following lemma gives the parametrization of the lattice points on $\gamma_z$ and two geometric identities. Part~\ref{sl:b} gives absolute convergence of line sums in Lemma~\ref{lem:convergence}. Part~\ref{sl:d} is the main tool in Section~\ref{sec:dyadic}.
\begin{figure}[t]
\centering
\begin{tikzpicture}[scale=0.9]
  \draw[gray!55,thin] (-1.5,-1.5) grid (4.5,4.5);
  \draw[gray!75,->] (-1.7,0) -- (4.7,0) node[right] {$y_1$};
  \draw[gray!75,->] (0,-1.7) -- (0,4.7) node[above] {$y_2$};
  \draw[red!75!black,thick] (3.25,-1.5) -- (0.25,4.5);
  \draw[blue!70!black,thick,->] (0,0) -- (2,1) node[midway,below right,font=\small] {$z$};
  \draw[blue!70!black,thick,->] (2,1) -- (1,3) node[midway,left,font=\small] {$v_z$};
  \foreach \p in {(2,1),(1,3),(3,-1)} {\fill \p circle (0.09);}
  \node[font=\small,anchor=west] at (2.1,1.3) {$z=(2,1)$};
  \node[font=\small,anchor=west] at (1.1,3.2) {$z+v_z$};
  \node[font=\small,anchor=west] at (3.1,-0.8) {$z-v_z$};
  \node[red!75!black,font=\small,anchor=south west] at (0.3,4.1) {$\gamma_z$};
\end{tikzpicture}
\caption{The selected line $\gamma_z$ for $z=(2,1)$. It passes through $z$ and is perpendicular to $z$. Its lattice points are $z+kv_z$, $k\in\Z$, with $v_z=Jz=(-1,2)$. The origin is not on $\gamma_z$, and $z$ is its lattice point closest to the origin.}
\label{fig:selected}
\end{figure}
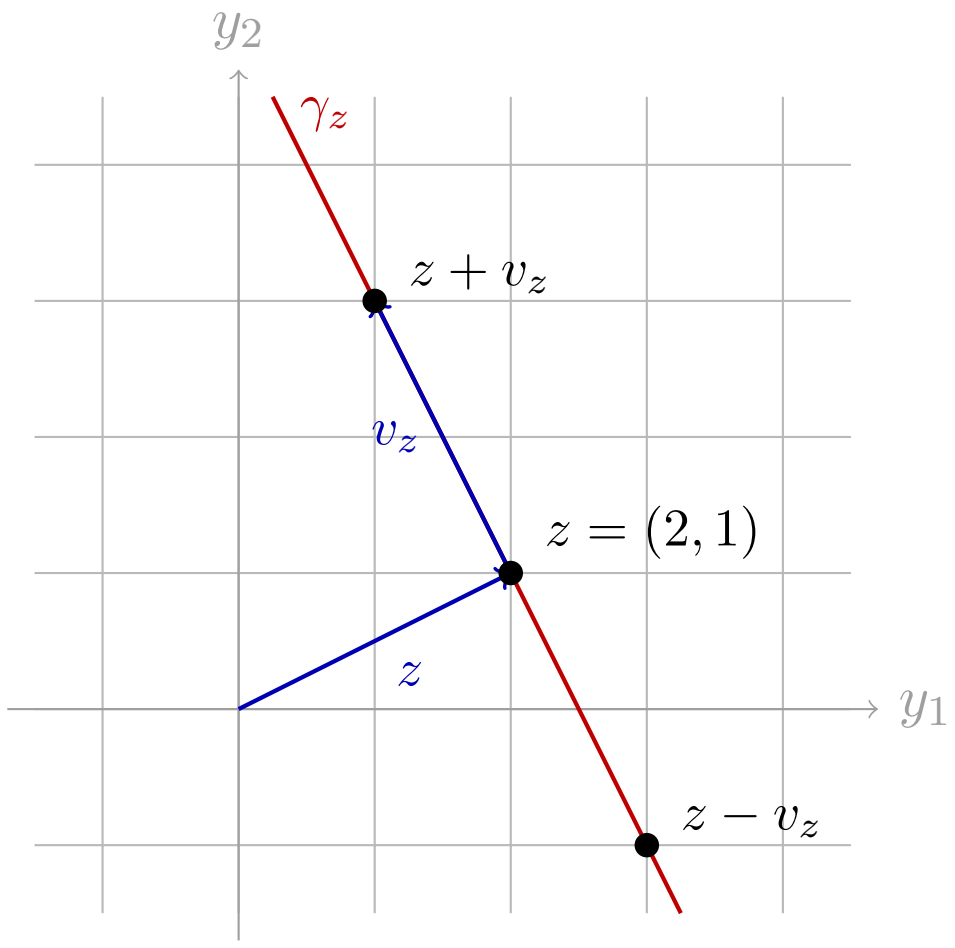
\begin{lemma}\label{lem:selected}
For every $z\in\Z^2$, the following hold.
\begin{enumerate}[label=\textup{(\roman*)},leftmargin=2.2em]
\item\label{sl:a}
We have $\gamma_z=z+\R v_z$, and $v_z$ is primitive.
The lattice points of $\gamma_z$ are given by
\[
\gamma_z\cap\Z^2=\{z+kv_z:k\in\Z\},
\]
where distinct integers $k$ give distinct points.

\item\label{sl:b}
For every $k\in\Z$, we have
\[
|z+kv_z|^2=|z|^2+k^2|v_z|^2,
\qquad |v_z|\ge1.
\]
In particular, $|z+kv_z|\ge|k|$.

\item\label{sl:c}
Every lattice point of $\gamma_z$ other than $z$
has norm strictly larger than $|z|$.
\end{enumerate}

For every $z\in\Z^2\setminus\{0\}$,
the following additional statements hold.
\begin{enumerate}[label=\textup{(\roman*)},leftmargin=2.2em,resume]
\item\label{sl:d}
A point $y\in\R^2$ belongs to $\gamma_z$
if and only if $|2z-y|^2=|y|^2$.

\item\label{sl:e}
We have $0\notin\gamma_z$.
\item\label{sl:f}
The point $z$ is the unique point of $\gamma_z$ closest to the origin.
\end{enumerate}
\end{lemma}
Figure~\ref{fig:selected} illustrates the selected line
and its primitive direction vector for $z=(2,1)$.
\begin{proof}
For part~\ref{sl:a}, first let $z\in\Z^2\setminus\{0\}$ and $y\in\R^2$.
By Lemma~\ref{lem:J}, the vectors $z$ and $Jz$ are nonzero and orthogonal. Therefore they form a basis of $\R^2$, so we can write $y-z=sz+tJz$ with $s,t\in\R$. Taking the dot product with $z$ gives
\[
z\cdot y-|z|^2=z\cdot(y-z)=s|z|^2.
\]
Since $|z|^2>0$, the equation $z\cdot y=|z|^2$ holds if and only if $s=0$.
This is equivalent to $y\in z+\R Jz$. Thus the two descriptions in \eqref{eq:gammaz} agree.

By Lemma~\ref{lem:bezout}\ref{bz:a}, the vector $z/g(z)$ is primitive.
Part~\ref{bz:d} of the same lemma then shows that $v_z=J(z/g(z))$ is primitive. Since $g(z)>0$, the vectors $v_z$ and $Jz$ have the same direction. Hence
\[
\gamma_z=z+\R Jz=z+\R v_z.
\]

For $z=0$, the definitions give $\gamma_0=\R\theta=0+\R v_0$. The vector $v_0=\theta$ is primitive by assumption. In both cases, Lemma~\ref{lem:rational}, applied with $p=z$ and $w=v_z$, gives
\[
\gamma_z\cap\Z^2=\{z+kv_z:k\in\Z\}.
\]
The same lemma shows that distinct integers $k$ give distinct points.

For part~\ref{sl:b}, let $z\in\Z^2$ and $k\in\Z$. We first note that $z\cdot v_z=0$. For $z\ne0$, this follows from $z\cdot Jz=0$ and $v_z=Jz/g(z)$. For $z=0$, the dot product vanishes directly. Consequently, we have
\[
|z+kv_z|^2
=|z|^2+2k\,z\cdot v_z+k^2|v_z|^2
=|z|^2+k^2|v_z|^2.
\]
Since $v_z$ is a nonzero integer vector, $|v_z|^2$ is a positive integer.
Thus $|v_z|^2\ge1$, and therefore $|v_z|\ge1$. The identity just proved now gives $|z+kv_z|^2\ge k^2$, so $|z+kv_z|\ge|k|$.

For part~\ref{sl:c}, let $y\in\gamma_z\cap\Z^2$ with $y\ne z$. By part~\ref{sl:a}, we have $y=z+kv_z$ for an integer $k\ne0$. Part~\ref{sl:b} gives
\[
|y|^2-|z|^2=k^2|v_z|^2>0.
\]
Hence $|y|>|z|$.

For part~\ref{sl:d}, let $z\ne0$ and $y\in\R^2$. Expanding the squared norm gives
\[
|2z-y|^2-|y|^2
=4|z|^2-4z\cdot y
=4(|z|^2-z\cdot y).
\]
Thus $|2z-y|^2=|y|^2$ holds if and only if $z\cdot y=|z|^2$. By \eqref{eq:gammaz}, this is exactly the condition $y\in\gamma_z$.

For part~\ref{sl:e}, let $z\ne0$. Then $z\cdot0=0\ne|z|^2$, so the origin does not satisfy the equation defining $\gamma_z$. 
Therefore $0\notin\gamma_z$.

For part~\ref{sl:f}, let $y\in\gamma_z$. The defining equation gives $z\cdot(y-z)=0$, so $|y|^2=|z|^2+|y-z|^2$. Hence $|y|\ge|z|$, with equality if and only if $y=z$.
\end{proof}
\begin{corollary}\label{cor:bijection}
The map $z\mapsto\gamma_z$ is a bijection from $\Z^2$ onto $\Tstar$.
\end{corollary}
\begin{proof}
Surjectivity holds by the definition \eqref{eq:Tstar}. For injectivity, let $z,z'\in\Z^2$ satisfy $\gamma_z=\gamma_{z'}$. If both $z$ and $z'$ are nonzero, then by Lemma~\ref{lem:selected}\ref{sl:f} each of them is the unique point of the common line closest to the origin, so $z=z'$. If one of them is zero, say $z=0$, then $0\in\gamma_0=\gamma_{z'}$, and Lemma~\ref{lem:selected}\ref{sl:e} forces $z'=0$.
\end{proof}
For $z\in\Z^2\setminus\{0\}$, the map $y\mapsto2z-y$ on $\R^2$ is reflection about $z$. It maps $\gamma_z$ onto itself. By Lemma~\ref{lem:selected}\ref{sl:d}, the equality $|2z-y|=|y|$ holds precisely for $y\in\gamma_z$.
\subsection{Absolute convergence of line sums}
Lemma~\ref{lem:selected}\ref{sl:a}-\ref{sl:b} gives a uniform bound for the selected line sums under polynomial decay.
\begin{lemma}\label{lem:convergence}
Let $s>1$, and let $f\colon\Z^2\to\C$ satisfy
\[
|f(z)|\le C_f(1+|z|)^{-s}
\qquad(z\in\Z^2)
\]
for some finite constant $C_f\ge0$. The constant may depend on $f$ and $s$, but not on $z$. Then, for every $z\in\Z^2$, we have
\begin{equation}\label{eq:linebound}
\sum_{y\in\gamma_z\cap\Z^2}|f(y)|
\le C_f\Bigl(1+\frac{2}{s-1}\Bigr).
\end{equation}
In particular, every selected line sum converges absolutely.

Let $\sigma(z):=(\gamma_z\cap\Z^2)\setminus\{z\}$. Then
\begin{equation}\label{eq:outward}
f(z)=Pf(\gamma_z)-\sum_{y\in\sigma(z)}f(y)
\qquad(z\in\Z^2).
\end{equation}
\end{lemma}
\begin{proof}
Fix $z\in\Z^2$. By Lemma~\ref{lem:selected}\ref{sl:a}, every lattice point $y$ on $\gamma_z$ can be written as $y=z+kv_z$ for exactly one integer $k$.
Part~\ref{sl:b} gives $|z+kv_z|\ge|k|$.
The decay assumption therefore gives
\[
\begin{aligned}
\sum_{y\in\gamma_z\cap\Z^2}|f(y)|
&=\sum_{k\in\Z}|f(z+kv_z)|\le C_f\Bigl(1+2\sum_{k\ge1}(1+k)^{-s}\Bigr).
\end{aligned}
\]
Since $t\mapsto t^{-s}$ is decreasing on $[1,\infty)$, we have
\[
\sum_{k\ge1}(1+k)^{-s}
\le\int_1^\infty t^{-s}\,dt
=\frac{1}{s-1}.
\]
This proves \eqref{eq:linebound}.

The term $k=0$ corresponds to the point $z$. Absolute convergence allows us to separate this term, giving
\[
Pf(\gamma_z)=f(z)+\sum_{y\in\sigma(z)}f(y).
\]
Rearranging this proves \eqref{eq:outward}.
\end{proof}
Every $f\in X_\varepsilon$ satisfies the bound in Lemma~\ref{lem:convergence} with $s=1+\varepsilon>1$, by \eqref{eq:Xeps}.
Thus every selected line sum converges absolutely.

Rapidly decreasing functions and functions satisfying $f(z)=o(|z|^{-2})$ satisfy
\[
|f(z)|\le C_f(1+|z|)^{-2}
\]
for all sufficiently large $|z|$. The lattice points inside the remaining bounded disk form a finite set. By choosing $C_f$ at least as large as
$(1+|z|)^2|f(z)|$ at each of these points, the estimate holds for every $z\in\Z^2$. Lemma~\ref{lem:convergence}, with $s=2$, therefore gives absolute convergence of the selected line sums.
\subsection{Divisibility of sums of two squares}

The counterexample construction and the dilation argument use two elementary divisibility properties of sums of squares. We recall these for the sake of readability.

For integers $a,b$ and a positive integer $n$, we write $a\equiv b\pmod n$ if $n\mid(a-b)$. Equivalently, $a$ and $b$ have the same remainder after division by $n$. 
\begin{lemma}\label{lem:squares}
Let $a,b\in\Z$.
\begin{enumerate}[label=\textup{(\roman*)},leftmargin=2.2em]
\item\label{sq:c}
If $4$ divides $a^2+b^2$, then $a$ and $b$
are both even.

\item\label{sq:d}
If $3$ divides $a^2+b^2$, then $3$ divides
both $a$ and $b$.
\end{enumerate}
\end{lemma}
\begin{proof}
For part~\ref{sq:c}, every integer is of the form $2k$ or $2k+1$ for some $k\in\Z$. The identities
\[
(2k)^2=4k^2,
\qquad
(2k+1)^2=4(k^2+k)+1
\]
show that an integer square is congruent to $0$ or $1$ modulo $4$. Thus $4\mid a^2+b^2$ forces both squares to be $0$ modulo $4$, so $a$ and $b$ are even.

For part~\ref{sq:d}, every integer is of the form $3k$, $3k+1$, or $3k-1$ for some $k\in\Z$. We have
\[
(3k)^2=9k^2,
\qquad
(3k\pm1)^2=3(3k^2\pm2k)+1.
\]
Hence an integer square is congruent to $0$ or $1$ modulo $3$.
Thus $3\mid a^2+b^2$ forces $3\mid a$ and $3\mid b$.
\end{proof}
\section{The counterexample}\label{sec:dyadic}
We first construct the values of $f_\theta$ at nonzero lattice points.
For $z\in\Z^2\setminus\{0\}$, Lemma~\ref{lem:J} gives
\begin{equation}\label{eq:twopoints}
z\cdot(z\pm Jz)=|z|^2,
\qquad
|z\pm Jz|^2=|z|^2+|Jz|^2=2|z|^2.
\end{equation}
Thus $z+Jz$ and $z-Jz$ lie on $\gamma_z$ and have squared norm twice that of $z$. This suggests choosing the support away from the origin on circles whose squared radii are powers of $2$.
\subsection{Lattice points on circles}

\begin{definition}\label{def:shells}
For an integer $m\ge0$ let
\[
D_m:=\{\,z\in\Z^2 : |z|^2=2^m\,\},\qquad D:=\bigcup_{m\ge0}D_m .
\]
Geometrically, $D_m$ is the intersection of $\Z^2$ with the circle of radius $2^{m/2}$ centered at the origin. The sets $D_m$ are pairwise disjoint, and $0\notin D$.
\end{definition}
To determine how the selected lines intersect $D$, we first identify all lattice points in each $D_m$.
\begin{lemma}\label{lem:shells}
For every integer $n\ge0$, we have
\[
D_{2n}=\{(\pm2^n,0),\,(0,\pm2^n)\},\qquad D_{2n+1}=\{(\pm2^n,\pm2^n)\}.
\]
For every integer $m\ge0$ and every $y\in D_m$, we also have $D_m=\{y,-y,Jy,-Jy\}$. Moreover, $D\subseteq\Sfour$.
\end{lemma}
The two signs in the formula for $D_{2n+1}$ are chosen independently.
\begin{proof}
Solving $a^2+b^2=1$ and $a^2+b^2=2$ in integers gives
\[
D_0=\{(\pm1,0),(0,\pm1)\},
\qquad
D_1=\{(\pm1,\pm1)\}.
\]
Let $m\ge2$ and $(a,b)\in D_m$. Since $4$ divides $a^2+b^2=2^m$, Lemma~\ref{lem:squares}\ref{sq:c} shows that $a$ and $b$ are even. Dividing by $4$ gives $(a/2)^2+(b/2)^2=2^{m-2}$, so $(a/2,b/2)\in D_{m-2}$.
Conversely, if $y\in D_{m-2}$, then $|2y|^2=4|y|^2=2^m$, so $2y\in D_m$. Hence
\[
D_m=2D_{m-2}
:=\{2y:y\in D_{m-2}\}
\qquad(m\ge2).
\]
Applying this identity $n$ times, we obtain
\[
D_{2n}=2^nD_0,\qquad D_{2n+1}=2^nD_1.
\]
Substituting the expressions for $D_0$ and $D_1$ gives the stated formulas.

Now fix $m\ge0$ and $y\in D_m$. By Lemma~\ref{lem:J}, the points $y,-y,Jy,-Jy$ belong to $D_m$. They are distinct because $y$ and $Jy$ are nonzero and orthogonal. Since $D_m$ has four points, $D_m=\{y,-y,Jy,-Jy\}$. The explicit formulas also show that $D\subseteq\Sfour$.
\end{proof}

\subsection{Intersections of selected lines with \texorpdfstring{$D$}{D}}

By \eqref{eq:twopoints}, a selected line with base point $z\in D_m$ contains the two points $z\pm Jz\in D_{m+1}$. To compute line sums we must know that there are no further intersections with $D$, and we must also determine whether selected lines with base points outside $D$ intersect $D$. We combine Lemma~\ref{lem:selected}\ref{sl:d} with Lemma~\ref{lem:shells}.
\begin{lemma}\label{lem:incidence}
Let $z\in\Z^2\setminus\{0\}$.
\begin{enumerate}[label=\textup{(\roman*)},leftmargin=2.2em]
\item\label{inc:a} If $z\notin D$, then $\gamma_z\cap D=\emptyset$.
\item\label{inc:b} If $z\in D_m$ for an integer $m\ge0$, then
\[
\gamma_z\cap D=\{z,z+Jz,z-Jz\}.
\]
These points are distinct, and
$z\pm Jz\in D_{m+1}$.
\end{enumerate}
\end{lemma}
\begin{proof}
Let $y\in\gamma_z\cap D$, and choose an integer
$m'\ge0$ with $y\in D_{m'}$. The point $2z-y$ has integer coordinates, and Lemma~\ref{lem:selected}\ref{sl:d} gives $|2z-y|^2=|y|^2=2^{m'}$. Hence $2z-y\in D_{m'}$, and Lemma~\ref{lem:shells} gives
\[
2z-y\in\{y,-y,Jy,-Jy\}.
\]
We consider the following cases. The possibility $2z-y=-y$ contradicts $z\ne0$. If $2z-y=y$, then $z=y\in D$. In the remaining cases,
\[
2z=y+\eta Jy,
\qquad \eta\in\{1,-1\}.
\]
That is, $2z=(I+\eta J)y$. Since $J^2=-I$ and $\eta^2=1$, we have $(I-\eta J)(I+\eta J)=2I$. Multiplying by $I-\eta J$ gives $2y=2z-2\eta Jz$, so $y=z-\eta Jz\in\{z+Jz,z-Jz\}$. Since $y$ and $Jy$ are orthogonal and have the same norm, the identity $2z=y+\eta Jy$ gives
\[
4|z|^2=|y|^2+|Jy|^2=2^{m'+1}.
\]
Hence, $|z|^2=2^{m'-1}$. If $m'=0$, then $|z|^2=1/2$, which is impossible because $z$ has integer coordinates. Thus $m'\ge1$, and the definition of $D_{m'-1}$ gives $z\in D_{m'-1}\subseteq D$.

We have proved that $\gamma_z\cap D\ne\emptyset$ implies $z\in D$.
Therefore, if $z\notin D$, then the intersection is empty, which proves part~\ref{inc:a}. The same argument also gives
\[
\gamma_z\cap D\subseteq\{z,z+Jz,z-Jz\}
\qquad(z\ne0).
\]

Now let $z\in D_m$. Then $z\in\gamma_z\cap D$, and \eqref{eq:twopoints} gives $z\pm Jz\in\gamma_z\cap D_{m+1}$. Together with the inclusion just proved, this gives part~\ref{inc:b}. The three points are distinct because $Jz\ne0$.
\end{proof}
\subsection{Construction}
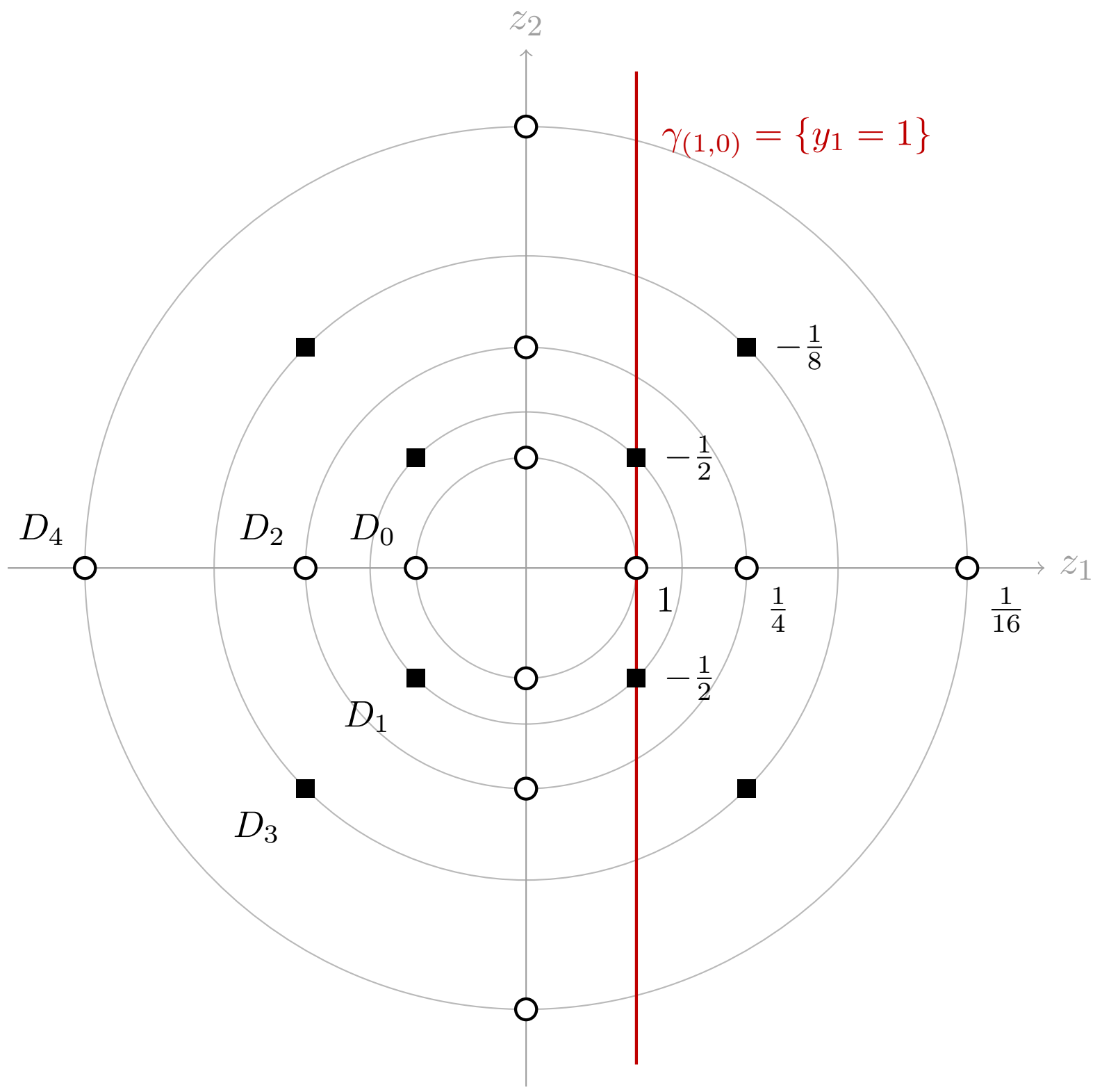
\begin{figure}[t]
\centering
\begin{tikzpicture}[scale=1.05]
  \foreach \r in {1,1.41421,2,2.82843,4} {\draw[gray!55,thin] (0,0) circle (\r);}
  \draw[gray!75,->] (-4.7,0) -- (4.7,0) node[right] {$z_1$};
  \draw[gray!75,->] (0,-4.7) -- (0,4.7) node[above] {$z_2$};
  \draw[red!75!black,thick] (1,-4.5) -- (1,4.5);
  \node[red!75!black,font=\small,anchor=west] at (1.1,3.9) {$\gamma_{(1,0)}=\{y_1=1\}$};
  \foreach \p in {(1,0),(-1,0),(0,1),(0,-1),(2,0),(-2,0),(0,2),(0,-2),(4,0),(-4,0),(0,4),(0,-4)} {\draw[fill=white,thick] \p circle (0.095);}
  \foreach \p in {(1,1),(1,-1),(-1,1),(-1,-1),(2,2),(2,-2),(-2,2),(-2,-2)} {\fill \p +(-0.085,-0.085) rectangle +(0.085,0.085);}
  \node[font=\small,anchor=north west] at (1.05,-0.05) {$1$};
  \node[font=\small,anchor=west] at (1.12,1) {$-\tfrac12$};
  \node[font=\small,anchor=west] at (1.12,-1) {$-\tfrac12$};
  \node[font=\small,anchor=north west] at (2.05,-0.05) {$\tfrac14$};
  \node[font=\small,anchor=west] at (2.12,2) {$-\tfrac18$};
  \node[font=\small,anchor=north west] at (4.05,-0.05) {$\tfrac1{16}$};
  \node[font=\small,anchor=south east] at (-1.05,0.08) {$D_0$};
  \node[font=\small,anchor=north east] at (-1.1,-1.08) {$D_1$};
  \node[font=\small,anchor=south east] at (-2.05,0.08) {$D_2$};
  \node[font=\small,anchor=north east] at (-2.1,-2.08) {$D_3$};
  \node[font=\small,anchor=south east] at (-4.05,0.08) {$D_4$};
\end{tikzpicture}
\caption{The sets $D_0,\ldots,D_4$. The numerical labels give values of $h$ at selected points.}
\label{fig:shells}
\end{figure}

Consider a function that vanishes outside $D$ and takes a constant value $a_m$ on each $D_m$. By Lemma~\ref{lem:incidence}\ref{inc:b}, its line sum over $\gamma_z$, with $z\in D_m$, equals $a_m+2a_{m+1}$. For these line sums to vanish, we require
\begin{equation}\label{eq:am-recurrence}
a_{m+1}=-\frac12a_m
\qquad(m\ge0).
\end{equation}
Iterating \eqref{eq:am-recurrence} with $a_0=1$ gives $a_m=(-1/2)^m$ for every $m\ge0$.
\begin{definition}\label{def:h}
For every integer $m\ge0$, put $a_m:=(-1/2)^m=(-1)^m2^{-m}$. Define $h\colon\Z^2\to\R$ by
\[
h(z):=
\begin{cases}
a_m, & z\in D_m,\\
0, & z\notin D.
\end{cases}
\]
Since $\theta$ is primitive, $\gamma_0=\R\theta$ is a coordinate axis if and only if $\theta\in\{(\pm1,0),(0,\pm1)\}$, and it is one of the lines $z_2=\pm z_1$ if and only if $\theta\in\{(\pm1,\pm1)\}$. In every other case $\gamma_0$ meets $\Sfour$ only at the origin. Define
\[
c_\theta:=
\begin{cases}
-\dfrac83, & \theta\in\{(\pm1,0),(0,\pm1)\},\\[2mm]
\dfrac43, & \theta\in\{(\pm1,\pm1)\},\\[2mm]
0, & \text{otherwise}.
\end{cases}
\]
Finally, put
\[
f_\theta:=h+c_\theta\one_{\{0\}},
\]
where $\one_{\{0\}}$ equals $1$ at the origin and $0$ elsewhere.
\end{definition}
Figure~\ref{fig:shells} shows $D_0,\ldots,D_4$
and some values of $h$.

The following lemma computes the line sums of $h$ and justifies the choice of $c_\theta$.
\begin{lemma}\label{lem:hprops}
The function $h$ has the following properties.
\begin{enumerate}[label=\textup{(\roman*)},leftmargin=2.2em]
\item\label{hp:a}
$h(0)=0$, $h(1,0)=1$, and
$\supp h=D\subseteq\Sfour$.

\item\label{hp:b}
For $z\in D$, we have $|h(z)|=|z|^{-2}$.
Consequently, $|h(z)|\le|z|^{-2}$
for every $z\ne0$.

\item\label{hp:c}
We have $\|h\|_{\ell^1}=8$.

\item\label{hp:d}
We have $Ph(\gamma_z)=0$
for every $z\in\Z^2\setminus\{0\}$.

\item\label{hp:e}
We have $Ph(\gamma_0)=-c_\theta$.
\end{enumerate}
\end{lemma}
\begin{proof}
For part~\ref{hp:a}, every $a_m$ is nonzero, and $h$ vanishes outside $D$.
Hence $\supp h=D$. Since $0\notin D$, we have $h(0)=0$. Also, $(1,0)\in D_0$ gives $h(1,0)=a_0=1$. The inclusion $D\subseteq\Sfour$ follows
from Lemma~\ref{lem:shells}.  

For part~\ref{hp:b}, let $z\in D_m$. Then $|z|^2=2^m$, so
\[
|h(z)|=|a_m|=2^{-m}=|z|^{-2}.
\]
If $z\notin D$ and $z\ne0$, then $|h(z)|=0\le|z|^{-2}$.

For part~\ref{hp:c}, the sets $D_m$ are pairwise disjoint and each contains four points by Lemma~\ref{lem:shells}. Hence
\[
\|h\|_{\ell^1}
=4\sum_{m\ge0}|a_m|
=4\sum_{m\ge0}2^{-m}
=8.
\]

For part~\ref{hp:d}, let $z\in\Z^2\setminus\{0\}$. If $z\notin D$, then
$\gamma_z\cap D=\emptyset$ by Lemma~\ref{lem:incidence}\ref{inc:a}. Since $h$ vanishes outside $D$, we obtain $Ph(\gamma_z)=0$. If $z\in D_m$, then
Lemma~\ref{lem:incidence}\ref{inc:b} and \eqref{eq:am-recurrence} give
\[
Ph(\gamma_z)=a_m+2a_{m+1}=0.
\]

For part~\ref{hp:e}, suppose first that $\gamma_0$ is neither a coordinate axis nor one of the lines $z_2=\pm z_1$. The line $\gamma_0$ meets each coordinate axis and each of the lines $z_2=\pm z_1$ only at the origin.
Since $D\subseteq\Sfour$ and $0\notin D$, we have $\gamma_0\cap D=\emptyset$. Thus $Ph(\gamma_0)=0$.

If $\gamma_0$ is a coordinate axis, Lemma~\ref{lem:shells} shows that,
for every integer $n\ge0$, it contains exactly two points of $D_{2n}$ and no points of $D_{2n+1}$. Hence
\[
Ph(\gamma_0)
=2\sum_{n\ge0}a_{2n}
=2\sum_{n\ge0}4^{-n}
=\frac{2}{1-\frac14}
=\frac83.
\]

If $\gamma_0$ is one of the lines $z_2=\pm z_1$, Lemma~\ref{lem:shells} shows that, for every integer $n\ge0$, it contains exactly two points of $D_{2n+1}$ and no points of $D_{2n}$. Hence
\[
Ph(\gamma_0)
=2\sum_{n\ge0}a_{2n+1}
=-\sum_{n\ge0}4^{-n}
=-\frac{1}{1-\frac14}
=-\frac43.
\]
By Definition~\ref{def:h}, these values give $Ph(\gamma_0)=-c_\theta$ in all cases.
\end{proof}
\subsection{Proof of Theorem~A\ref{A:a}}
\begin{proof}
Recall from Definition~\ref{def:h} that $f_\theta=h+c_\theta\one_{\{0\}}$.
By Lemma~\ref{lem:hprops}\ref{hp:a}, we have $f_\theta(1,0)=1$, so $f_\theta$ is nonzero. Since $\supp h=D$ and $f_\theta$ differs from $h$ only at the origin, we have
\[
\supp f_\theta\subseteq D\cup\{0\}\subseteq\Sfour.
\]
Since $h(0)=0$, Lemma~\ref{lem:hprops}\ref{hp:c} gives
\[
\|f_\theta\|_{\ell^1}
=\|h\|_{\ell^1}+|c_\theta|
=8+|c_\theta|<\infty.
\]
Thus every line sum of $f_\theta$ converges absolutely. For $z\ne0$, Lemma~\ref{lem:hprops}\ref{hp:b} also gives
\[
|f_\theta(z)|=|h(z)|\le|z|^{-2}.
\]

For $z\ne0$, the origin does not belong to $\gamma_z$ by Lemma~\ref{lem:selected}\ref{sl:e}. Thus $f_\theta(y)=h(y)$ for every
$y\in\gamma_z\cap\Z^2$. Lemma~\ref{lem:hprops}\ref{hp:d} therefore gives
\[
Pf_\theta(\gamma_z)=Ph(\gamma_z)=0.
\]
Since $\gamma_0$ contains the origin, Lemma~\ref{lem:hprops}\ref{hp:e} gives
\[
Pf_\theta(\gamma_0)=Ph(\gamma_0)+c_\theta=0.
\]
Therefore $Pf_\theta(\gamma)=0$ for every $\gamma\in\Tstar$.

Let $0<\varepsilon\le1$. For every nonzero lattice point $z$, we have $|z|\ge1$, so
\[
(1+|z|)^{1+\varepsilon}|f_\theta(z)|
\le (1+|z|)^2|z|^{-2}
=\Bigl(1+\frac1{|z|}\Bigr)^2
\le4.
\]
At the origin, we have $|f_\theta(0)|=|c_\theta|\le8/3$. Thus the supremum in \eqref{eq:Xeps} is finite, and $f_\theta\in X_\varepsilon$.

Now let $\varepsilon>1$ and put $z_n:=(2^n,0)$ for integers $n\ge0$. By Lemma~\ref{lem:shells}, we have $z_n\in D_{2n}$. Since $z_n\ne0$, Definition~\ref{def:h} gives
\[
|f_\theta(z_n)|=|h(z_n)|=|a_{2n}|=4^{-n}.
\]
Consequently,
\[
\begin{aligned}
(1+|z_n|)^{1+\varepsilon}|f_\theta(z_n)|
&=(1+2^n)^{1+\varepsilon}4^{-n}\ge2^{n(\varepsilon-1)}
\longrightarrow\infty
\qquad(n\to\infty).
\end{aligned}
\]
Thus $f_\theta\notin X_\varepsilon$.
\end{proof}
\section{An infinite-dimensional null space}\label{sec:dilation}
By Lemma~\ref{lem:hprops}\ref{hp:d}, the function $h$ has zero sums on every selected line with nonzero base point. Its sums on the coordinate axes and the lines $z_2=\pm z_1$ are nonzero, as computed in
Lemma~\ref{lem:hprops}\ref{hp:e}. By taking differences of dilates of $h$, we obtain functions whose line sums vanish on every line through the origin.

We use the dilation factor $3$ for two reasons. First,  dilation by $3$ preserves the zero sums of $h$ on all selected lines
with nonzero base point as shown in Lemma~\ref{lem:dilation}. Second, the dilated supports $3^jD$, $j\ge0$, are also pairwise disjoint. We use this disjointness in Proposition~\ref{prop:Gj} to prove that the constructed functions are linearly independent.
\begin{definition}\label{def:T3}
For every integer $j\ge0$, we write
\[
3^j\Z^2:=\{3^jy:y\in\Z^2\}.
\]
Thus $3\Z^2$ consists of the lattice points whose two coordinates are divisible by $3$.

For a function $u\colon\Z^2\to\C$, we define $T_3u\colon\Z^2\to\C$ by
\begin{equation}\label{eq:T3}
(T_3u)(z):=
\begin{cases}
u(z/3), & z\in3\Z^2,\\
0, & z\notin3\Z^2.
\end{cases}
\end{equation}
\end{definition}
We define the iterates by $T_3^0u=u$ and
$T_3^{j+1}u=T_3(T_3^ju)$ for integers $j\ge0$.
For every $j\ge0$ and $z\in\Z^2$, we have
\[
(T_3^ju)(z)=
\begin{cases}
u(z/3^j), & z\in3^j\Z^2,\\
0, & z\notin3^j\Z^2.
\end{cases}
\]

We prove this formula by induction on $j$.
The case $j=0$ follows from $T_3^0u=u$.
Suppose that the formula holds for some $j\ge0$.
If $z\notin3\Z^2$, then \eqref{eq:T3} gives
$(T_3^{j+1}u)(z)=0$.
If $z\in3\Z^2$, then
\[
(T_3^{j+1}u)(z)=(T_3^ju)(z/3).
\]
We apply the induction hypothesis at $z/3$.
Since $z/3\in3^j\Z^2$ if and only if
$z\in3^{j+1}\Z^2$, the resulting value is
$u(z/3^{j+1})$ when $z\in3^{j+1}\Z^2$
and zero otherwise.
This proves the formula for $j+1$.

\begin{lemma}\label{lem:dilation}
Let $u\in\ellone$. Then $T_3u\in\ellone$ and
\[
\|T_3u\|_{\ell^1}=\|u\|_{\ell^1}.
\]
Moreover, the following hold.
\begin{enumerate}[label=\textup{(\roman*)},leftmargin=2.2em]
\item\label{dl:a}
If $z\in\Z^2\setminus3\Z^2$, then
$\gamma_z\cap3\Z^2=\emptyset$.
Consequently, $P(T_3u)(\gamma_z)=0$.

\item\label{dl:b}
For every $z'\in\Z^2\setminus\{0\}$, we have
\[
P(T_3u)(\gamma_{3z'})=Pu(\gamma_{z'}).
\]

\item\label{dl:c}
For every line $L\subset\R^2$ through the origin,
we have
\[
P(T_3u)(L)=Pu(L).
\]
\end{enumerate}
\end{lemma}
\begin{proof}
Since $T_3u$ vanishes outside $3\Z^2$ and $y\mapsto3y$ is a bijection from $\Z^2$ onto $3\Z^2$, we have
\[
\|T_3u\|_{\ell^1}
=\sum_{x\in3\Z^2}|u(x/3)|
=\sum_{y\in\Z^2}|u(y)|
=\|u\|_{\ell^1}<\infty.
\]
In particular, all line sums of $u$ and $T_3u$ converge absolutely.

For part~\ref{dl:a}, let $z\in\Z^2\setminus3\Z^2$. Suppose that $x\in\gamma_z\cap3\Z^2$, and write $x=3y$ with $y\in\Z^2$. The equation defining $\gamma_z$ gives
\[
z_1^2+z_2^2=|z|^2=z\cdot x=3(z\cdot y).
\]
Thus $3$ divides $z_1^2+z_2^2$. By Lemma~\ref{lem:squares}\ref{sq:d},
both coordinates of $z$ are divisible by $3$. This contradicts $z\notin3\Z^2$, so $\gamma_z\cap3\Z^2=\emptyset$. Since $T_3u$ vanishes outside $3\Z^2$, we obtain $P(T_3u)(\gamma_z)=0$.

For part~\ref{dl:b}, let $z'\in\Z^2\setminus\{0\}$. For $y\in\Z^2$, the equation $(3z')\cdot(3y)=|3z'|^2$ is equivalent to $z'\cdot y=|z'|^2$.
Thus $3y\in\gamma_{3z'}$ if and only if $y\in\gamma_{z'}$.
Hence $y\mapsto3y$ is a bijection from $\gamma_{z'}\cap\Z^2$ onto
$\gamma_{3z'}\cap3\Z^2$. Since $T_3u$ vanishes outside $3\Z^2$,
reindexing the absolutely convergent sum gives
\[
\begin{aligned}
P(T_3u)(\gamma_{3z'})
&=\sum_{x\in\gamma_{3z'}\cap3\Z^2}u(x/3)=\sum_{y\in\gamma_{z'}\cap\Z^2}u(y)
=Pu(\gamma_{z'}).
\end{aligned}
\]

For part~\ref{dl:c}, let $L$ be a line through the origin.
Since $L$ is a linear subspace of $\R^2$,
we have $3y\in L$ if and only if $y\in L$.
Using \eqref{eq:T3} and reindexing the absolutely
convergent sum, we obtain
\[
\begin{aligned}
P(T_3u)(L)
&=\sum_{x\in L\cap3\Z^2}u(x/3)=\sum_{y\in L\cap\Z^2}u(y)
=Pu(L).
\end{aligned}
\]
\end{proof}
We now take differences of consecutive dilates of $h$
to prove Theorem~A\ref{A:b}.
\begin{proposition}\label{prop:Gj}
For every integer $j\ge0$, put
\[
h_j:=T_3^jh,
\qquad
G_j:=h_j-h_{j+1}.
\]
The functions $G_j$ are real-valued, belong to $X_1\cap\ellone$, and are linearly independent. For every $j\ge0$, we have
\[
G_j(0)=0,
\qquad
\|G_j\|_{\ell^1}=16,
\qquad
|G_j(z)|\le9^{j+1}|z|^{-2}
\quad(z\ne0).
\]
Moreover, each $G_j$ satisfies
\[
PG_j(\gamma_z)=0
\qquad(z\in\Z^2\setminus\{0\}),
\]
and $PG_j(L)=0$ for every line $L\subset\R^2$ through the origin. 
\end{proposition}
\begin{proof}
As $h$ is real-valued, $G_j$ is also real-valued. Since $h_j(0)=h(0)=0$, we have $G_j(0)=0$. Repeated application of Lemma~\ref{lem:dilation} gives
\[
\|h_j\|_{\ell^1}=\|h\|_{\ell^1}=8.
\]

We prove by induction on $j$ that
\begin{equation}\label{eq:hjdata}
Ph_j(\gamma_z)=0
\qquad(z\in\Z^2\setminus\{0\}).
\end{equation}
For $j=0$, we have $h_0=h$, so \eqref{eq:hjdata} follows from Lemma~\ref{lem:hprops}\ref{hp:d}. Suppose that \eqref{eq:hjdata} holds for some $j\ge0$,
and fix $z\in\Z^2\setminus\{0\}$.
If $z\notin3\Z^2$, then
Lemma~\ref{lem:dilation}\ref{dl:a} gives
\[
Ph_{j+1}(\gamma_z)
=P(T_3h_j)(\gamma_z)=0.
\]
If $z\in3\Z^2$, put $z':=z/3$.
Then $z'\in\Z^2\setminus\{0\}$.
By Lemma~\ref{lem:dilation}\ref{dl:b}
and the induction hypothesis, we obtain
\[
Ph_{j+1}(\gamma_z)
=Ph_j(\gamma_{z'})=0.
\]
Thus \eqref{eq:hjdata} holds for $j+1$, and hence for every $j\ge0$. Since $G_j=h_j-h_{j+1}$, we conclude that
\[
PG_j(\gamma_z)
=Ph_j(\gamma_z)-Ph_{j+1}(\gamma_z)=0
\qquad(z\ne0).
\]

For every line $L$ through the origin, repeated application of
Lemma~\ref{lem:dilation}\ref{dl:c} gives $Ph_j(L)=Ph_{j+1}(L)=Ph(L)$.
Therefore $PG_j(L)=0$.

Since $\supp h=D$, we have $\supp h_j=3^jD$ for every $j\ge0$. We show that these supports are pairwise disjoint. Let $0\le j<k$ and suppose that
$z\in3^jD\cap3^kD$. Then
\[
|z|^2=9^j2^m=9^k2^n
\]
for some integers $m,n\ge0$. Dividing by $9^j$, we obtain
\[
2^m=9^{k-j}2^n.
\]
Thus $3$ divides $2^m$, which is impossible because $2^m\equiv(-1)^m\not\equiv0\pmod3$. Therefore the supports $3^jD$, $j\ge0$,
are pairwise disjoint. In particular, $h_j$ and $h_{j+1}$ have disjoint
supports, so
\[
\|G_j\|_{\ell^1}
=\|h_j\|_{\ell^1}+\|h_{j+1}\|_{\ell^1}
=16.
\]
For $z\in3^jD$, Lemma~\ref{lem:hprops}\ref{hp:b} gives
\[
|h_j(z)|
=|h(z/3^j)|
=9^j|z|^{-2}.
\]
Outside $3^jD$, the function $h_j$ vanishes. Since $h_j$ and $h_{j+1}$ have disjoint supports, at most one of $h_j(z)$ and $h_{j+1}(z)$ is nonzero. Hence
\[
|G_j(z)|\le9^{j+1}|z|^{-2}
\qquad(z\ne0).
\]
Every nonzero lattice point satisfies $|z|\ge1$, so
\[
(1+|z|)^2|G_j(z)|
\le9^{j+1}\Bigl(1+\frac1{|z|}\Bigr)^2
\le4\cdot9^{j+1}.
\]
Together with $G_j(0)=0$, this proves $G_j\in X_1$ by \eqref{eq:Xeps}.

To prove linear independence, suppose that 
\[
\sum_{j=0}^{M}b_jG_j=0
\]
holds with $b_j\in\C$ and $b_M\ne0$.
Put $q:=(3^{M+1},0)$.
Since $q=3^{M+1}(1,0)$, we have
\[
h_{M+1}(q)=h(1,0)=1.
\]
Thus $q\in\supp h_{M+1}=3^{M+1}D$. The pairwise disjointness of the supports gives $h_j(q)=0$ for $0\le j\le M$.
Consequently,
\[
G_j(q)=0\quad(0\le j<M),
\qquad
G_M(q)=-1.
\]
Evaluating the relation at $q$ gives $-b_M=0$,
contradicting $b_M\ne0$.
Therefore the functions $G_j$ are linearly independent.
\end{proof}
\section{Uniqueness for functions supported on four lines}
\label{sec:rigidity}
The counterexample $f_\theta$ satisfies $|z|^2|f_\theta(z)|=1$ on $D$.
We now prove exterior uniqueness for functions supported on $\Sfour$ under the stronger condition $f(z)=o(|z|^{-2})$. We first determine the intersections $\gamma_z\cap\Sfour$.
\begin{lemma}\label{lem:fourlines}
The following hold.
\begin{enumerate}[label=\textup{(\roman*)},leftmargin=2.2em]
\item\label{fl:a} Let $z\in\Sfour\setminus\{0\}$. Then we have
\[
\gamma_z\cap\Sfour=\{z,\,z+Jz,\,z-Jz\}.
\]
The three points are distinct. Moreover, $|z\pm Jz|=\sqrt2\,|z|$.
\item\label{fl:b} Let $z\in\Z^2\setminus\{0\}$. Then $\gamma_z\cap\Sfour$ has at most four points.
\end{enumerate}
\end{lemma}
\begin{proof}
The set $\Sfour$ consists of the lattice points on the four lines
\begin{align*}
    \ell_1&=\{y\in \R^2:y_2=0\},\qquad\ell_2=\{y\in \R^2:y_1=0\},\\\ell_3&=\{y\in \R^2:y_2=y_1\},\qquad\ell_4=\{y\in \R^2:y_2=-y_1\}.
\end{align*}

For part~\ref{fl:a}, let $z\in\Sfour\setminus\{0\}$. The rotation $J$ interchanges $\ell_1$ and $\ell_2$, and interchanges $\ell_3$ and $\ell_4$. It also maps $\Z^2$ onto itself. Thus $J\Sfour=\Sfour$. We also have $J\gamma_z=\gamma_{Jz}$. Indeed, for $y\in\R^2$, Lemma~\ref{lem:J} gives $Jz\cdot Jy=z\cdot y$ and $|Jz|^2=|z|^2$, so $y\in\gamma_z$ if and only if $Jy\in\gamma_{Jz}$ by \eqref{eq:gammaz}. Since
\[
J(z\pm Jz)=Jz\pm J(Jz),
\]
the intersection formula in part~\ref{fl:a} holds for $Jz$ whenever it holds for $z$. Therefore it suffices to consider $z=(a,0)$ and $z=(a,a)$, where $a\in\Z\setminus\{0\}$.

First, let $z=(a,0)$. By \eqref{eq:gammaz}, 
\[
\gamma_z=\{y\in\R^2:y_1=a\}.
\]
This line does not meet $\ell_2$. Its intersections with $\ell_1$, $\ell_3$, and $\ell_4$ give
\[
\gamma_z\cap\Sfour
=\{(a,0),(a,a),(a,-a)\}
=\{z,z+Jz,z-Jz\},
\]
since $Jz=(0,a)$.

Next, let $z=(a,a)$. In this case,
\[
\gamma_z=\{y\in\R^2:y_1+y_2=2a\}.
\]
This line does not meet $\ell_4$, since $a\ne0$. Its intersections with $\ell_3$, $\ell_1$, and $\ell_2$ give
\[
\gamma_z\cap\Sfour
=\{(a,a),(2a,0),(0,2a)\}
=\{z,z-Jz,z+Jz\},
\]
since $Jz=(-a,a)$.

Applying $J$ also covers the other coordinate axis and the other diagonal. Thus the intersection formula in part~\ref{fl:a}
holds for every $z\in\Sfour\setminus\{0\}$. The three points are distinct because $Jz\ne0$. Finally, \eqref{eq:twopoints} gives
$|z\pm Jz|=\sqrt2\,|z|$.

For part~\ref{fl:b}, let $z\in\Z^2\setminus\{0\}$. Each $\ell_i$ passes through the origin, whereas $\gamma_z$ does not, by
Lemma~\ref{lem:selected}\ref{sl:e}. Thus $\gamma_z$ differs from every $\ell_i$. Two distinct lines meet in at most one point, so $\gamma_z\cap\Sfour$ has at most four points.
\end{proof}
If $\supp f\subseteq\Sfour$, then Lemma~\ref{lem:fourlines}\ref{fl:b} shows that $Pf(\gamma_z)$ is a finite sum whenever $z\ne0$. Thus no decay assumption is needed for these line sums.
\begin{theorem}[Exterior uniqueness on four lines]
\label{thm:rigidity}
Let $\alpha\ge0$ and let $f\colon\Z^2\to\C$ satisfy $ \supp f\subseteq\Sfour$. Suppose that $Pf(\gamma_z)=0$ for every $z\in(\Sfour\cap E_\alpha)\setminus\{0\}$. Then, for every $z\in(\Sfour\cap E_\alpha)\setminus\{0\}$, we have
\begin{equation}\label{eq:recurrence}
f(z)=2f(2z)+f(2Jz)+f(-2Jz).
\end{equation}
Moreover, for every such $z$ and every integer $n\ge0$, we have
\begin{equation}\label{eq:iterate}
|z|^2|f(z)|
\le
\max_{0\le j\le3}
|2^nJ^jz|^2\,|f(2^nJ^jz)|.
\end{equation}

If, in addition, $f(x)=o(|x|^{-2})$ as $|x|\to\infty$, then $Pf(\gamma_0)$ converges absolutely and $f=0$ on $E_\alpha\setminus\{0\}$.
Under this decay assumption, if $\alpha=0$ and $Pf(\gamma_0)=0$, then also $f(0)=0$.
\end{theorem}
\begin{proof}
Let $w\in(\Sfour\cap E_\alpha)\setminus\{0\}$. The hypothesis gives $Pf(\gamma_w)=0$. Since $f$ vanishes outside $\Sfour$, Lemma~\ref{lem:fourlines}\ref{fl:a} gives
\begin{equation}\label{eq:atw}
0=Pf(\gamma_w)=f(w)+f(w+Jw)+f(w-Jw).
\end{equation}
Put $u:=w+Jw$ and $v:=w-Jw$. By Lemma~\ref{lem:fourlines}\ref{fl:a} we have $u,v\in\Sfour$ and \[
|u|=|v|=\sqrt2\,|w|>|w|\ge\alpha.
\] Thus $u,v\in(\Sfour\cap E_\alpha)\setminus\{0\}$. Hence the hypothesis also gives $Pf(\gamma_u)=Pf(\gamma_v)=0$.

Using $J^2=-I$, we obtain
\begin{align*}
u+Ju=2Jw, \quad u-Ju=2w,\quad
v+Jv=2w,  \quad v-Jv=-2Jw.
\end{align*}
By Lemma~\ref{lem:fourlines}\ref{fl:a} and the equalities $Pf(\gamma_u)=Pf(\gamma_v)=0$, we have
\begin{equation}\label{eq:atuv}
f(u)+f(2Jw)+f(2w)=0,\quad
f(v)+f(2w)+f(-2Jw)=0.
\end{equation}
Substituting the expressions for $f(u)$ and $f(v)$ from \eqref{eq:atuv} into \eqref{eq:atw}, we obtain 
\begin{align*}
    f(w)&=-f(u)-f(v)=2f(2w)+f(2Jw)+f(-2Jw).
\end{align*}
Since $w$ was arbitrary, this proves \eqref{eq:recurrence}.

Fix $z\in(\Sfour\cap E_\alpha)\setminus\{0\}$. For every integer $n\ge0$, define
\[
M_n:=\max_{0\le j\le3}|f(2^nJ^jz)|.
\]
Both $J$ and multiplication by $2$ map $\Sfour$ into itself.
Hence every point $2^nJ^jz$, with $j\in\{0,1,2,3\}$, belongs to $\Sfour$
and has norm
\[
|2^nJ^jz|=2^n|z|\ge|z|\ge\alpha.
\]
These points are nonzero. Applying \eqref{eq:recurrence} and using $-J=J^3$, we obtain
\[
\begin{aligned}
|f(2^nJ^jz)|
&\le 2|f(2^{n+1}J^jz)|
   +|f(2^{n+1}J^{j+1}z)|+|f(2^{n+1}J^{j+3}z)|\le4M_{n+1}.
\end{aligned}
\]
Since $J^4=I$, every integer power of $J$ equals one of $I,J,J^2,J^3$. Each absolute value is therefore bounded by $M_{n+1}$. Taking the maximum over $j$ gives $M_n\le4M_{n+1}$. Iterating this inequality gives
$M_0\le4^nM_n$ for every integer $n\ge0$. 

Since $|f(z)|\le M_0$ and $|2^nJ^jz|^2=4^n|z|^2$, we obtain
\[
\begin{aligned}
|z|^2|f(z)|
&\le4^n|z|^2M_n=\max_{0\le j\le3}
  |2^nJ^jz|^2\,|f(2^nJ^jz)|.
\end{aligned}
\]
This proves \eqref{eq:iterate}.

Assume that $f(x)=o(|x|^{-2})$ as $|x|\to\infty$.
For every $\eta>0$, there is $R>0$ such that
\[
|x|^2|f(x)|\le\eta
\qquad(x\in\Z^2,\ |x|\ge R).
\] Fix $z\in(\Sfour\cap E_\alpha)\setminus\{0\}$
and $\eta>0$. Choose an integer $n\ge0$ such that $2^n|z|\ge R$. Each point $2^nJ^jz$, with $j\in\{0,1,2,3\}$, has norm $2^n|z|\ge R$.
Hence every term in the maximum in \eqref{eq:iterate} is at most $\eta$, and we obtain
\[
|z|^2|f(z)|\le\eta.
\]
Since $\eta>0$ is arbitrary and $|z|>0$, we conclude that $f(z)=0$.
Together with $\supp f\subseteq\Sfour$, this proves $f=0$ on $E_\alpha\setminus\{0\}$.

Since $f(x)=o(|x|^{-2})$ as $|x|\to\infty$, there is a finite constant $C_f\ge0$ such that
\[
|f(x)|\le C_f(1+|x|)^{-2}
\qquad(x\in\Z^2),
\]
as explained after Lemma~\ref{lem:convergence}. That lemma, with $s=2$, shows that $Pf(\gamma_0)$ converges absolutely.
If $\alpha=0$, we have already proved that
$f=0$ on $E_0\setminus\{0\} =\Z^2\setminus\{0\}$. Since $\gamma_0$ contains the origin, we consequently have
\[
Pf(\gamma_0)=f(0).
\]
If, in addition, $Pf(\gamma_0)=0$, then $f(0)=0$.
\end{proof}
\begin{proof}[Proof of Theorem~A\ref{A:c}]
Theorem~\ref{thm:rigidity} gives $f=0$ on $E_\alpha\setminus\{0\}$. If $0\in E_\alpha$, then $\alpha=0$ and the hypothesis also gives $Pf(\gamma_0)=0$. Applying Theorem~\ref{thm:rigidity} with $\alpha=0$
and $Pf(\gamma_0)=0$, we also obtain $f(0)=0$. Hence $f$ vanishes at every point of $E_\alpha$.
\end{proof}

\subsection*{Acknowledgments}
The work of S.R.J.\ and J.R.\ was supported by the Research Council of Finland through the Flagship of Advanced Mathematics for Sensing, Imaging and Modelling (decision number 359183) and by the Jenny and Antti Wihuri Foundation. In addition, J.R.\ was supported by the Emil Aaltonen Foundation.

\subsection*{Data availability statement}
Data sharing is not applicable to this article, as no datasets were generated or analyzed during the current study.

\subsection*{Conflict of interest}
The authors declare no potential conflicts of interest with respect to the research, authorship, or publication of this article.

\subsection*{Declaration of generative AI use}
During the preparation of this manuscript, the authors
used ChatGPT Astra 6 and Claude Opus 5.5/Fable 5.1 for assistance with mathematical exposition, discussion of proofs and mathematical ideas, and LaTeX editing. The authors reviewed the AI-assisted content.
\bibliography{math} 

@article {Novikov:Sharma:2026,
    AUTHOR = {Novikov, Roman G. and Sharma, Basant Lal},
     TITLE = {On discrete {X}-ray transform},
   JOURNAL = {Bull. Sci. Math.},
  FJOURNAL = {Bulletin des Sciences Math\'ematiques},
    VOLUME = {213},
      YEAR = {2026},
     PAGES = {Paper No. 103940},
      ISSN = {0007-4497,1952-4773},
   MRCLASS = {44A12 (46F12 65N21 65N22 65R10)},
  MRNUMBER = {5134829},
       DOI = {10.1016/j.bulsci.2026.103940},
       URL = {https://doi.org/10.1016/j.bulsci.2026.103940},
}

@article {Strichartz:1982,
    AUTHOR = {Strichartz, Robert S.},
     TITLE = {Radon inversion---variations on a theme},
   JOURNAL = {Amer. Math. Monthly},
  FJOURNAL = {American Mathematical Monthly},
    VOLUME = {89},
      YEAR = {1982},
    NUMBER = {6},
     PAGES = {377--384, 420--423},
      ISSN = {0002-9890,1930-0972},
   MRCLASS = {44A15 (28A35 92A07)},
  MRNUMBER = {660917},
MRREVIEWER = {Ulrich\ Oppel},
       DOI = {10.2307/2321649},
       URL = {https://doi.org/10.2307/2321649},
}

@article {Gupta:Novikov:2026,
    AUTHOR = {Gupta, Pranav and Novikov, Roman G.},
     TITLE = {Discrete non-abelian {X}-ray transforms},
   JOURNAL = {J. Geom. Anal.},
  FJOURNAL = {Journal of Geometric Analysis},
    VOLUME = {36},
      YEAR = {2026},
    NUMBER = {1},
     PAGES = {Paper No. 37, 21},
      ISSN = {1050-6926,1559-002X},
   MRCLASS = {44A12 (65N21 65N22 65R10)},
  MRNUMBER = {5000395},
MRREVIEWER = {Lalit\ Mohan\ Upadhyaya},
       DOI = {10.1007/s12220-025-02281-0},
       URL = {https://doi.org/10.1007/s12220-025-02281-0},
}

@book {Jones:Jones:book,
    AUTHOR = {Jones, Gareth A. and Jones, J. Mary},
     TITLE = {Elementary number theory},
    SERIES = {Springer Undergraduate Mathematics Series},
 PUBLISHER = {Springer-Verlag London, Ltd., London},
      YEAR = {1998},
     PAGES = {xiv+301},
      ISBN = {3-540-76197-7},
   MRCLASS = {11-01},
  MRNUMBER = {1610533},
MRREVIEWER = {Charles\ J.\ Parry},
       DOI = {10.1007/978-1-4471-0613-5},
       URL = {https://doi.org/10.1007/978-1-4471-0613-5},
}

@article {Novikov:Xu:2026,
    AUTHOR = {Novikov, Roman G. and Xu, Tianli},
     TITLE = {On non-uniqueness of phase retrieval for discrete and
              continuous {Fourier} transforms},
   JOURNAL = {J. Inverse Ill-Posed Probl.},
  FJOURNAL = {Journal of Inverse and Ill-posed Problems},
  YEAR    = {2026},
   DOI = {10.1515/jiip-2026-0119},
       URL = {https://doi.org/10.1515/jiip-2026-0119},
}

@article {Abouelaz:2011:tori,
  AUTHOR = {Abouelaz, Ahmed and Rouvi{\`e}re, Fran{\c{c}}ois},
     TITLE = {Radon transform on the torus},
   JOURNAL = {Mediterr. J. Math.},
  FJOURNAL = {Mediterranean Journal of Mathematics},
    VOLUME = {8},
      YEAR = {2011},
    NUMBER = {4},
     PAGES = {463--471},
      ISSN = {1660-5446,1660-5454},
   MRCLASS = {53C65 (44A12)},
  MRNUMBER = {2860679},
MRREVIEWER = {Valeri\u i\ Vladimirovich\ Volchkov},
       DOI = {10.1007/s00009-010-0111-7},
       URL = {https://doi.org/10.1007/s00009-010-0111-7},
}

@unpublished{Rim:2026,
  author  = {Rim, Donsub},
  title   = {Explicit oblique inverses of the approximate
             discrete {Radon} transform},
  year    = {2026},
  note    = {Preprint, arXiv:2610.23363},
  doi     = {10.48550/arXiv.2610.23363},
  url     = {https://arxiv.org/abs/2610.23363}
}

@article {Li:Ren:Rim:2023,
    AUTHOR = {Li, Weilin and Ren, Kui and Rim, Donsub},
     TITLE = {A range characterization of the single-quadrant {ADRT}},
   JOURNAL = {Math. Comp.},
  FJOURNAL = {Mathematics of Computation},
    VOLUME = {92},
      YEAR = {2023},
    NUMBER = {339},
     PAGES = {283--306},
      ISSN = {0025-5718,1088-6842},
   MRCLASS = {65R10 (44A12 92C55)},
  MRNUMBER = {4496966},
MRREVIEWER = {Fatma\ Terzioglu},
       DOI = {10.1090/mcom/3750},
       URL = {https://doi.org/10.1090/mcom/3750},
}

@article {Rim:2020,
    AUTHOR = {Rim, Donsub},
     TITLE = {Exact and fast inversion of the approximate discrete {R}adon
              transform from partial data},
   JOURNAL = {Appl. Math. Lett.},
  FJOURNAL = {Applied Mathematics Letters. An International Journal of Rapid
              Publication},
    VOLUME = {102},
      YEAR = {2020},
     PAGES = {106159, 5},
      ISSN = {0893-9659,1873-5452},
   MRCLASS = {65T50 (44A12)},
  MRNUMBER = {4039172},
       DOI = {10.1016/j.aml.2019.106159},
       URL = {https://doi.org/10.1016/j.aml.2019.106159},
}

@article {Press:2006,
    AUTHOR = {Press, William H.},
     TITLE = {Discrete {R}adon transform has an exact, fast inverse and
              generalizes to operations other than sums along lines},
   JOURNAL = {Proc. Natl. Acad. Sci. USA},
  FJOURNAL = {Proceedings of the National Academy of Sciences of the United
              States of America},
    VOLUME = {103},
      YEAR = {2006},
    NUMBER = {51},
     PAGES = {19249--19254},
      ISSN = {0027-8424,1091-6490},
   MRCLASS = {65R10 (44A12 94A08)},
  MRNUMBER = {2270446},
       DOI = {10.1073/pnas.0609228103},
       URL = {https://doi.org/10.1073/pnas.0609228103},
}

@article {Brady:1998,
    AUTHOR = {Brady, Martin L.},
     TITLE = {A fast discrete approximation algorithm for the {R}adon
              transform},
   JOURNAL = {SIAM J. Comput.},
  FJOURNAL = {SIAM Journal on Computing},
    VOLUME = {27},
      YEAR = {1998},
    NUMBER = {1},
     PAGES = {107--119},
      ISSN = {0097-5397,1095-7111},
   MRCLASS = {44A12 (44-04 65R10 68U10)},
  MRNUMBER = {1614880},
MRREVIEWER = {Gabor\ T.\ Herman},
       DOI = {10.1137/S0097539793256673},
       URL = {https://doi.org/10.1137/S0097539793256673},
}

@unpublished {Salo:Meaney:Koskela:Railo:2026,
    AUTHOR = {Salo, Ella and Meaney, Alexander and Koskela, Olli
              and Railo, Jesse},
     TITLE = {Torus computed tomography for experimental data},
      YEAR = {2026},
      NOTE = {Preprint, arXiv:2608.16453},
       URL = {https://arxiv.org/abs/2608.16453},
}

@article{Boman:1993,
  author  = {Boman, Jan},
  title   = {An example of non-uniqueness for a generalized
             {Radon} transform},
   JOURNAL = {J. Anal. Math.},
  FJOURNAL = {Journal d'Analyse Math\'ematique},
  volume  = {61},
  year    = {1993},
  pages   = {395--401},
  doi     = {10.1007/BF02788850}
}

@article {Goncharov:Novikov:2019,
    AUTHOR = {Goncharov, F. O. and Novikov, R. G.},
     TITLE = {A breakdown of injectivity for weighted ray transforms in
              multidimensions},
   JOURNAL = {Ark. Mat.},
  FJOURNAL = {Arkiv f\"or Matematik},
    VOLUME = {57},
      YEAR = {2019},
    NUMBER = {2},
     PAGES = {333--371},
      ISSN = {0004-2080,1871-2487},
   MRCLASS = {44A12 (53C65 65R32)},
  MRNUMBER = {4018757},
MRREVIEWER = {Fritz\ Keinert},
       DOI = {10.4310/ARKIV.2019.v57.n2.a5},
       URL = {https://doi.org/10.4310/ARKIV.2019.v57.n2.a5},
}

@article {Railo:2020,
    AUTHOR = {Railo, Jesse},
     TITLE = {Fourier analysis of periodic {R}adon transforms},
   JOURNAL = {J. Fourier Anal. Appl.},
  FJOURNAL = {The Journal of Fourier Analysis and Applications},
    VOLUME = {26},
      YEAR = {2020},
    NUMBER = {4},
     PAGES = {Paper No. 64, 27},
      ISSN = {1069-5869,1531-5851},
   MRCLASS = {44A12 (42B05 45Q05 46F12)},
  MRNUMBER = {4128058},
MRREVIEWER = {S.\ K.\ Upadhyay},
       DOI = {10.1007/s00041-020-09775-1},
       URL = {https://doi.org/10.1007/s00041-020-09775-1},
}

@article {Armitage:1994,
    AUTHOR = {Armitage, D. H.},
     TITLE = {A non-constant continuous function on the plane whose integral
              on every line is zero},
   JOURNAL = {Amer. Math. Monthly},
  FJOURNAL = {American Mathematical Monthly},
    VOLUME = {101},
      YEAR = {1994},
    NUMBER = {9},
     PAGES = {892--894},
      ISSN = {0002-9890,1930-0972},
   MRCLASS = {30E99 (26B35)},
  MRNUMBER = {1300495},
MRREVIEWER = {J.\ Kr\'al},
       DOI = {10.2307/2975138},
       URL = {https://doi.org/10.2307/2975138},
}

@article {Zalcman:1982,
    AUTHOR = {Zalcman, Lawrence},
     TITLE = {Uniqueness and nonuniqueness for the {R}adon transform},
   JOURNAL = {Bull. London Math. Soc.},
  FJOURNAL = {The Bulletin of the London Mathematical Society},
    VOLUME = {14},
      YEAR = {1982},
    NUMBER = {3},
     PAGES = {241--245},
      ISSN = {0024-6093,1469-2120},
   MRCLASS = {42B10 (44A05 92A07)},
  MRNUMBER = {656606},
MRREVIEWER = {R.\ B.\ Burckel},
       DOI = {10.1112/blms/14.3.241},
       URL = {https://doi.org/10.1112/blms/14.3.241},
}

@article {Jensen:2004,
    AUTHOR = {Jensen, Sine R.},
     TITLE = {Sufficient conditions for the inversion formula for the
              {$k$}-plane {R}adon transform in {$\bold R^n$}},
   JOURNAL = {Math. Scand.},
  FJOURNAL = {Mathematica Scandinavica},
    VOLUME = {94},
      YEAR = {2004},
    NUMBER = {2},
     PAGES = {207--226},
      ISSN = {0025-5521,1903-1807},
   MRCLASS = {44A12},
  MRNUMBER = {2053741},
MRREVIEWER = {Alexander\ Katsevich},
       DOI = {10.7146/math.scand.a-14440},
       URL = {https://doi.org/10.7146/math.scand.a-14440},
}

@article {Lehtonen:Railo:Salo:2018,
    AUTHOR = {Lehtonen, Jere and Railo, Jesse and Salo, Mikko},
     TITLE = {Tensor tomography on {C}artan-{H}adamard manifolds},
   JOURNAL = {Inverse Problems},
  FJOURNAL = {Inverse Problems. An International Journal on the Theory and
              Practice of Inverse Problems, Inverse Methods and Computerized
              Inversion of Data},
    VOLUME = {34},
      YEAR = {2018},
    NUMBER = {4},
     PAGES = {044004, 27},
      ISSN = {0266-5611,1361-6420},
   MRCLASS = {44A12},
  MRNUMBER = {3781677},
MRREVIEWER = {Aleksander\ Denisiuk},
       DOI = {10.1088/1361-6420/aaaf85},
       URL = {https://doi.org/10.1088/1361-6420/aaaf85},
}

@article {Ilmavirta:Koskela:Railo:2020,
    AUTHOR = {Ilmavirta, Joonas and Koskela, Olli and Railo, Jesse},
     TITLE = {Torus computed tomography},
   JOURNAL = {SIAM J. Appl. Math.},
  FJOURNAL = {SIAM Journal on Applied Mathematics},
    VOLUME = {80},
      YEAR = {2020},
    NUMBER = {4},
     PAGES = {1947--1976},
      ISSN = {0036-1399,1095-712X},
   MRCLASS = {65R10 (42B05 44A12 46F12 65R32)},
  MRNUMBER = {4140552},
MRREVIEWER = {Aleksander\ Denisiuk},
       DOI = {10.1137/19M1268070},
       URL = {https://doi.org/10.1137/19M1268070},
}

@article {Ilmavirta:2015,
    AUTHOR = {Ilmavirta, Joonas},
     TITLE = {On {R}adon transforms on tori},
   JOURNAL = {J. Fourier Anal. Appl.},
  FJOURNAL = {The Journal of Fourier Analysis and Applications},
    VOLUME = {21},
      YEAR = {2015},
    NUMBER = {2},
     PAGES = {370--382},
      ISSN = {1069-5869,1531-5851},
   MRCLASS = {46F12 (44A12)},
  MRNUMBER = {3319537},
MRREVIEWER = {Takashi\ Takiguchi},
       DOI = {10.1007/s00041-014-9374-x},
       URL = {https://doi.org/10.1007/s00041-014-9374-x},
}

@article {Fokas:Novikov:1991,
    AUTHOR = {Fokas, Athanassios S. and Novikov, Roman G.},
     TITLE = {Discrete analogues of {$\overline\partial$}-equation and of
              {R}adon transform},
   JOURNAL = {C. R. Acad. Sci. Paris S\'er. I Math.},
  FJOURNAL = {Comptes Rendus de l'Acad\'emie des Sciences. S\'erie I.
              Math\'ematique},
    VOLUME = {313},
      YEAR = {1991},
    NUMBER = {2},
     PAGES = {75--80},
      ISSN = {0764-4442},
   MRCLASS = {39A12},
  MRNUMBER = {1119912},
}

@book {GIP2D,
    AUTHOR = {Paternain, Gabriel P. and Salo, Mikko and Uhlmann, Gunther},
     TITLE = {Geometric inverse problems---with emphasis on two dimensions},
    SERIES = {Cambridge Studies in Advanced Mathematics},
    VOLUME = {204},
      NOTE = {With a foreword by Andr\'as Vasy},
 PUBLISHER = {Cambridge University Press, Cambridge},
      YEAR = {2023},
     PAGES = {xxiv+344},
      ISBN = {978-1-316-51087-2},
   MRCLASS = {35-02 (35R30 53C65 58J32 92C55)},
  MRNUMBER = {4520155},
}

@book {Sharafutdinov:book,
    AUTHOR = {Sharafutdinov, V. A.},
     TITLE = {Integral geometry of tensor fields},
    SERIES = {Inverse and Ill-posed Problems Series},
 PUBLISHER = {VSP, Utrecht},
      YEAR = {1994},
     PAGES = {271},
      ISBN = {90-6764-165-0},
   MRCLASS = {53C65 (44A12 73D50)},
  MRNUMBER = {1374572},
MRREVIEWER = {Eric\ Grinberg},
       DOI = {10.1515/9783110900095},
       URL = {https://doi.org/10.1515/9783110900095},
}

@book {Helgason:Book,
    AUTHOR = {Helgason, Sigurdur},
     TITLE = {The {R}adon transform},
    SERIES = {Progress in Mathematics},
    VOLUME = {5},
   EDITION = {Second},
 PUBLISHER = {Birkh\"auser Boston, Inc., Boston, MA},
      YEAR = {1999},
     PAGES = {xiv+188},
      ISBN = {0-8176-4109-2},
   MRCLASS = {44A12 (22E30 22E46 43A85 58J40)},
  MRNUMBER = {1723736},
MRREVIEWER = {B.\ S.\ Rubin},
       DOI = {10.1007/978-1-4757-1463-0},
       URL = {https://doi.org/10.1007/978-1-4757-1463-0},
}
	
	\bibliographystyle{alpha}
\end{document}